\documentclass[12pt]{amsart}

\newcommand{\Z}{\mathbb Z}
\newcommand{\R}{\mathbb R}
\newcommand{\N}{\mathbb N}
\newcommand{\G}{\mathcal G}
\newcommand{\M}{\mathcal M}
\renewcommand{\L}{\mathcal L}

\newtheorem{theorem}{Theorem}

\newtheorem{lemma}[theorem]{Lemma}
\newtheorem{corollary}[theorem]{Corollary}
\newtheorem{definition}[theorem]{Definition}

\allowdisplaybreaks[4]

\begin{document}

\title{Non-MSF wavelets for the Hardy space $H^2(\R)$}
\author{Biswaranjan Behera} 
\address{Statistics and Mathematics Unit, Indian Statistical Institute, 203, B. T. Road,
Calcutta, 700108, INDIA}
\email{\tt biswa\_v@isical.ac.in, br\_behera@yahoo.com}
\date{May 16, 2002 \\ \hspace*{3.5mm} {2000 {\it Mathematics Subject Classification.} 42C40}}


\begin{abstract}
We prove three results on wavelets for the Hardy space $H^2(\R)$. All wavelets constructed so far for $H^2(\R)$ are MSF wavelets. We construct a family of $H^2$-wavelets which are not MSF. An equivalence relation on $H^2$-wavelets is introduced and it is shown that the corresponding equivalence classes are non-empty. Finally, we construct a family of $H^2$-wavelets with Fourier transform discontinuous at the origin.
\end{abstract}

\maketitle
\section{Introduction}

The classical Hardy space $H^2(\R)$ is the collection of all square integrable functions whose Fourier transform is supported in $\R^+=(0,\infty)$:
\[
H^2(\R):= \{f\in L^2(\R): \hat f(\xi)=0 \mbox{~for a.e. } \xi\leq 0 \},
\]
where $\hat f$ is the Fourier transform of $f$ defined by
\[
\hat f(\xi)=\int_{\R}f(x) e^{-i\xi x} dx.
\]
Clearly, $H^2(\R)$ is a closed subspace of $L^2(\R)$. A function $\psi\in H^2(\R)$ is said to be a wavelet for $H^2(\R)$ if the system of functions 
$\{\psi_{j,k}=2^{j/2}\psi(2^j\cdot-k):j,k\in\Z\}$ 
forms an orthonormal basis for $H^2(\R)$. We shall call such a $\psi$ an $H^2$-wavelet.

Two basic equations characterize all $H^2$-wavelets. The proof of the following theorem can be obtained from the corresponding result for the usual case of $L^2(\R)$ 
(see Theorem~6.4, Chapter~7 in~\cite{HW}).

\begin{theorem} 
\label{T.one}
A function $\psi\in H^2(\R)$, with $\|\psi\|_2=1$, is an $H^2$-wavelet if and only if 
\[
\sum_{j\in \Z} |\hat\psi(2^j\xi)|^2 = \chi_{\R^+}
\quad{\rm for~a.e.}~\xi\in\R 
\]
and
\[ 
\sum_{j\geq 0} \hat\psi(2^j\xi) \overline{\hat\psi(2^j(\xi +2q\pi))}= 0
\quad{\rm for~a.e.}~\xi\in\R~{\rm and~for~all}~q\in 2\Z+1.
\]
\end{theorem}

From the Paley-Wiener theorem it follows that there is no compactly supported function in
$H^2(\R)$ apart from the zero function, hence, there is no compactly supported $H^2$-wavelet. On the other hand, there exist $H^2$-wavelets with compactly supported Fourier transform. One such example is given by $\hat\psi=\chi_{[2\pi,4\pi]}$ which is the analogue of the Shannon wavelet for $L^2(\R)$. P.\ Auscher~\cite{AUS} proved that there is no $H^2$-wavelet satisfying the following regularity condition: $|\hat\psi|$ is continuous on $\R$ and $|\hat\psi(\xi)| = O((1+|\xi|)^{-\alpha-\frac{1}{2}})$ 
at $\infty$, for some $\alpha>0$. In particular, $H^2(\R)$ does not have a wavelet $\psi$ with $|\hat\psi|$ continuous and $\hat\psi$ is compactly supported. 

Analogous to the $L^2$ case, an $H^2$-wavelet $\psi$ will be called a {\it minimally supported frequency} (MSF) wavelet if $|\hat\psi|=\chi_K$ for some $K\subset\R^+$. Such wavelets were called {\it unimodular wavelets} in~\cite{HKLS} and {\it s-elementary wavelets} in~\cite{DL}. The associated set $K$ will be called an {\it $H^2$-wavelet set}. In this situation the set $K$ has Lebesgue measure $2\pi$.

There is a simple characterization of $H^2$-wavelet sets analogous to the $L^2$ case.
\begin{theorem} 
\label{T.two}
A set $K\subset\R^+$ is an $H^2$-wavelet set if and only if the following two conditions hold:
\begin{enumerate}
	\item[(i)] $\{K+2k\pi:k\in\Z\}$ is a partition of $\R$.
	\item[(ii)] $\{2^j K:j\in\Z\}$ is a partition of $\R^+$.
\end{enumerate}
\end{theorem}

In ~\cite{HKLS}, the authors proved that the only $H^2$-wavelet set which is an interval is $[2\pi, 4\pi]$. They also characterized all $H^2$-wavelet sets consisting of two disjoint intervals. In~\cite{B} (see also~\cite{ABM}) we proved a result on the structure of $H^2$-wavelet sets consisting of a finite number of intervals and, as an application,  characterized 3-interval $H^2$-wavelet sets. All these wavelet sets depend on a finite number of integral parameters which proves that there are countably many $H^2$-wavelet sets which are union of at most three disjoint intervals. We also constructed a family of 4-interval $H^2$-wavelet sets with some of the endpoints depending on a continuous real parameter, thereby proving the uncountability of such sets (see~\cite{ABM}). In the proof of Theorem~\ref{T.five} below, we exibit a family of $H^2$-wavelet sets with some of the endpoints depending on two independent continuous real parameters. Some more $H^2$-wavelet sets were constructed in~\cite{MAJ} where the author also proves the existence of an $H^2$-MSF wavelet $\psi$ such that $\psi\not\in L^p(\R)$ for $p<2$.

Up to the present time, all known wavelets for $H^2(\R)$ are MSF, i.e., the Fourier transform is the characteristic function of a subset of $\R^+$. In the next section, 
we construct a family of non-MSF $H^2$-wavelet. In section~3, we introduce an equivalence relation on the set of $H^2$-wavelets and explicitly construct examples of wavelets in each of the corresponding equivalence classes. In the last section, we construct a family of $H^2$-wavelets with Fourier transform discontinuous at the origin.
\section{The construction of non-MSF wavelets}

Our strategy of constructing the family of non-MSF wavelets of $H^2(\R)$ is the following. We start with an $H^2$-MSF wavelet so that $|\hat\psi|$ assumes the value 1 on its support. Then we add some more sets to the support of $\hat\psi$ and reassign values to $\hat\psi$ in such a manner that the equalities $\sum_{j\in\Z}|\hat\psi(2^j\xi)|^2=\chi_{\R^+}$ a.e. and 
$\sum_{k\in\Z}|\hat\psi(\xi+2k\pi)|^2=1$ a.e. are preserved, which are necessary conditions for $\psi$ to be an $H^2$-wavelet.

Fix $r\in\N$ and let $k$ be an integer satisfying $1\leq k< 2(2^r-1)$. Define the set
\[
K_{r,k} = 
\Bigl[\frac{2(k+1)}{2^{r+1}-1}\pi, \frac{2k}{2^r-1}\pi\Bigr]\cup
\Bigl[\frac{2^{r+1}k}{2^r-1}\pi, \frac{2^{r+2}(k+1)}{2^{r+1}-1}\pi\Bigr]
=A\cup B,~{\rm say}.
\]

Observe that the sets $A+2\cdot2^r\pi$ and $B$ are disjoint and their union is an interval of length $2\pi$ so that (i) in Theorem~\ref{T.two} is satisfied. Similarly, $2^r A$ and $B$ are disjoint and their union is the interval $[a,2a]$, where $a=\tfrac{2^{r+1}(k+1)}{2^{r+1}-1}\pi$, hence, (ii) in Theorem~\ref{T.two} is also satisfied. Therefore, $K_{r,k}$ is an $H^2$-wavelet set. In fact, 
$\{K_{r,k}:r\in\N, 1\leq k<2(2^r-1)\}$ is precisely the collection of all $H^2$-wavelet sets consisting of two disjoint intervals, as first observed in~\cite{HKLS}.

In particular, for $k=2^r-1$, we get the following family of $H^2$-wavelet sets:
\begin{equation}
\label{E.kr}
K_r=
\Bigl[\frac{2^{r+1}}{2^{r+1}-1}\pi, 2\pi\Bigr]\cup
\Bigl[2^{r+1}\pi, \frac{2^{2r+2}}{2^{r+1}-1}\pi\Bigr],
\quad r\in\N.
\end{equation}
Let us denote the intervals on the right hand side of~(\ref{E.kr}) by $I_r$ and $J_r$ respectively. Note that $\frac{2\pi}{3}\leq |I_r|<\pi$ and $\pi<|J_r|\leq\frac{4\pi}{3}$. We denote the Lebesgue measure of  a set $S$ by $|S|$. First of all, we observe that $2^{-1}I_r+2^{r+1}\pi\subset J_r$.

For $r\in\N$, define the function $\psi_r$ by
\begin{eqnarray}
\label{E.four}
\hat\psi_r(\xi)=
\begin{cases}
\frac{1}{\sqrt{2}} & {\rm if}~\xi\in I_r\cup (2^{-1}I_r)\cup(2^{-1}I_r+2^{r+1}\pi) \\
-\frac{1}{\sqrt 2} & {\rm if}~\xi\in I_r+2^{r+2}\pi \\
1 & {\rm if}~\xi\in J_r\setminus(2^{-1}I_r+2^{r+1}\pi) \\
0 & {\rm otherwise}.
\end{cases}
\end{eqnarray}

\begin{theorem}
\label{T.three}
For each $r\in\N$, $\psi_r$ is a wavelet for the Hardy space $H^2(\R)$.
\end{theorem}

Some preparation is needed before we prove Theorem~\ref{T.three}. Define the maps $\tau$ and $d$ as follows:
\[
\tau:\R\rightarrow [2\pi,4\pi], \quad \tau(x)=x+2k(x)\pi,
\]
\[
d:\R^+\rightarrow [2\pi,4\pi], \quad d(x)=2^{j(x)}x,
\]
where $k(x)$ and $j(x)$ are unique integers such that $x+2k(x)\pi$ and $2^{j(x)}x$ belong to $[2\pi,4\pi]$. 

We first prove the following lemma which gives useful information regarding the support of $\hat\psi_r$. This will be crucial for proving Theorem~\ref{T.three}.

\begin{lemma}
\label{L.one}
Let $E_r={\rm supp\ }\hat\psi_r=(2^{-1}I_r)\cup I_r\cup J_r\cup(I_r+2^{r+2}\pi)$.
\begin{enumerate} 
\item [(i)] If $\xi\in2^{-1}I_r$, 
  then $\xi+2k\pi\in E_r$ if and only if $k=0, 2^r$, 
  and $2^j\xi\in E_r$ if and only if $j=0,1$. 
\item [(ii)] If $\xi\in I_r$, 
  then $\xi+2k\pi\in E_r$ if and only if $k=0,2^{r+1}$, 
  and $2^j\xi\in E_r$ if and only if $j=0,-1$. 
\item [(iii)] If $\xi\in 2^{-1}I_r+2^{r+1}\pi$, 
  then $\xi+2k\pi\in E_r$ if and only if $k=0, -2^r$, 
  and $2^j\xi\in E_r$ if and only if $j=0,1$.
\item [(iv)] If $\xi\in J_r\setminus(2^{-1}I_r+2^{r+1}\pi)$, 
  then $\xi+2k\pi\in E_r$ if and only if $k=0$, 
  and $2^j\xi\in E_r$ if and only if $j=0$.
\item [(v)] If $\xi\in I_r+2^{r+2}\pi$, 
  then $\xi+2k\pi\in E_r$ if and only if $k=0,-2^{r+1}$, 
  and $2^j\xi\in E_r$ if and only if $j=0,-1$.
\end{enumerate}
\end{lemma}

\proof
Observe that $\tau(E)=\tau(E+2k\pi)$ and $d(F)=d(2^jF)$ for every $j,k\in\Z$ and for every $E\subset\R$, $F\subset\R^+$. Hence,
\begin{equation}
\label{E.tau1}
\tau(2^{-1}I_r+2^{r+1}\pi)=\tau(2^{-1}I_r),\ \tau(I_r+2^{r+2}\pi)=\tau(I_r),
\end{equation}
and 
\begin{equation}
\label{E.d1}
d(2^{-1}I_r)=d(I_r),\ d(2^{-1}I_r+2^{r+1}\pi)=d(I_r+2^{r+2}\pi).
\end{equation}
It also follows from the definition of the maps $\tau$ and $d$ that if $W$ is an $H^2$-wavelet set and $E,F\subset W$, then $\tau(E)\cap\tau(F)=\emptyset$ and 
$d(E)\cap d(F)=\emptyset$. Since $I_r\cup J_r$ is an $H^2$-wavelet set and $2^{-1}I_r+2^{r+1}\pi\subset J_r$, we have
\begin{equation}
\label{E.tau2}
\tau(I_r)\cap\tau(2^{-1}I_r+2^{r+1}\pi)=\emptyset,
\end{equation}
and
\begin{equation}
\label{E.tau3}
\tau(J_r\setminus(2^{-1}I_r+2^{r+1}\pi))\cap\tau(2^{-1}I_r+2^{r+1}\pi)=\emptyset.
\end{equation}

Form (\ref{E.tau1}), (\ref{E.tau2}) and (\ref{E.tau3}), we get 
\[
\tau(2^{-1}I_r)\cap\tau(I_r)=\emptyset,\
\tau(2^{-1}I_r)\cap\tau(I_r+2^{r+2}\pi)=\emptyset,
\]
and
\[
\tau(2^{-1}I_r)\cap\tau(J_r\setminus(2^{-1}I_r+2^{r+1}\pi))=\emptyset.
\]
Therefore, if $\xi\in2^{-1}I_r$, then $\xi+2k\pi\in E_r$ if and only if $k=0, 2^r$.

Similarly, we have
\begin{equation}
\label{E.d2}
d(I_r)\cap d(J_r)=\emptyset\quad {\rm and}\quad
d(I_r)\cap d(2^{-1}I_r+2^{r+1}\pi)=\emptyset.
\end{equation}
From (\ref{E.d1}) and (\ref{E.d2}) we get 
\[
d(2^{-1}I_r)\cap d(J_r)=\emptyset\quad {\rm and}\quad
d(2^{-1}I_r)\cap d(I_r+2^{r+2}\pi)=\emptyset.
\]
From this we obtain that if $\xi\in2^{-1}I_r$, then $2^j\xi\in E_r$ if and only if $j=0,1$.

We have proved (i) of the theorem. The proof of (ii)--(v) is similar.
\qed

\proof[Proof of Theorem~{\rm \ref{T.three}}]
In view of the characterization of $H^2$-wavelets (see Theorem~\ref{T.one}), we need to show the following:
\begin{enumerate}
\item[(a)] $\|\psi_r\|_2=1$. \\
\item[(b)] $\rho(\xi):=\sum_{j\in\Z}|\hat\psi_r(2^j\xi)|^2=
           \chi_{\R^+}$ for a.e. $\xi\in\R$. \\
\item[(c)] $t_q(\xi):= \sum_{j\geq 0}	  
           \hat\psi_r(2^j\xi)\overline{\hat\psi_r(2^j(\xi+2q\pi))}=0$ 
           for a.e. $\xi\in\R$ and for all $q\in2\Z+1$.
\end{enumerate}

{\it Proof of {\rm (a)}}. We have
\begin{eqnarray*}
\|\hat\psi_r\|_2^2 
& = & \int_\R |\hat\psi_r(\xi)|^2 d\xi \\
& = & \tfrac{1}{2}\Bigl(|I_r|+\tfrac{1}{2}|I_r|+\tfrac{1}{2}|I_r|+|I_r|\Bigr)
       + |J_r|-\tfrac{1}{2}|I_r|\\
& = & |I_r|+|J_r| \\
& = & 2\pi.
\end{eqnarray*}
Hence, $\|\psi\|_2=1$. 

{\it Proof of {\rm (b)}}. Observe that $\rho(\xi)=0$ if $\xi\leq 0$. Since $\rho(2\xi)=\rho(\xi)$ for a.e. $\xi$, it is enough to show that $\rho(\xi)=1$ on any set $E$ such that $d(E)=[2\pi,4\pi]$. $I_r\cup J_r$ is such a set since it is an $H^2$-wavelet set.

If $\xi\in I_r$, then by Lemma~\ref{L.one}(ii), $2^j\xi\in$ supp~$\hat\psi_r$ if and only if $j=-1,0$. Hence, 
$\rho(\xi)=|\hat\psi(\frac{\xi}{2})|^2+|\hat\psi(\xi)|^2
=(\frac{1}{\sqrt{2}})^2+(\frac{1}{\sqrt{2}})^2=1$. 

We write $J_r=(2^{-1}I_r+2^{r+1}\pi)\cup
\{J_r\setminus(2^{-1}I_r+2^{r+1}\pi)\} = M\cup L$, say. 
If  $\xi\in M$, then $2^j\xi\in$ supp~$\hat\psi_r$ if and only if $j=0,1$ 
(see Lemma~\ref{L.one}(iii)) so that 
$\rho(\xi)=|\hat\psi(\xi)|^2+|\hat\psi(2\xi)|^2
=(\frac{1}{\sqrt{2}})^2+(-\frac{1}{\sqrt{2}})^2=1$. For $\xi\in L$, no other dilate of $\xi$ is in the support of $\hat\psi$, hence, $\rho(\xi)=1$ a.e.

{\it Proof of {\rm (c)}}. Since $t_{-q}(\xi)=\overline{t_q(\xi-2q\pi)}$, it is enough to prove that $t_q=0$ a.e. for all positive and odd integer $q$. The term 
\[
\hat\psi(2^j\xi)\overline{\hat\psi(2^j(\xi+2q\pi))}
\]
is non-zero when both $2^j\xi$ and $2^j\xi+2\cdot 2^jq\pi$ are in the support of $\hat\psi$. Referring again to Lemma~\ref{L.one}, we observe that, this is possible if either $2^j q=2^r$ or $2^j q=2^{r+1}$. Since the integer $q$ is odd, either $j=r,q=1$ or $j=r+1,q=1$. In the first case, $2^j\xi\in2^{-1}I_r$ so that $2^j(\xi+2q\pi)\in2^{-1}I_r+2^{r+1}\pi$, $2^{j+1}\xi\in I_r$, and  
$2^{j+1}(\xi+2q\pi)\in I_r+2^{r+2}\pi$. 
Hence, 
$t_q(\xi)=(\frac{1}{\sqrt{2}})(\frac{1}{\sqrt{2}})
+(\frac{1}{\sqrt{2}})(-\frac{1}{\sqrt{2}})=0$. 
The second case is treated similarly. This completes the proof of the theorem.
\qed
\section{An equivalence relation}

In this section we shall introduce an equivalence relation on the collection of all wavelets of $H^2(\R)$ and show that each of the corresponding equivalence classes is non-empty.

Let $\psi$ be an $H^2$-wavelet. For $j\in\Z$, define the following closed subspaces of $H^2(\R)$:
$V_j =\overline{\rm span}\{\psi_{l,k}:\mbox{$l<j$},\mbox{$k\in\Z$}\}$. 
It is easy to verify that these subspaces satisfy the following properties:
\begin{enumerate}
\item[(i)] $V_j\subset V_{j+1}$ for all $j\in\Z$,
\item[(ii)] $f\in V_j$ if and only if $f(2\cdot)\in V_{j+1}$ for all $j\in\Z$,
\item[(iii)] $\cup_{j\in\Z} V_j$ is dense in $H^{2}(\R)$, $\cap_{j\in\Z}=\{0\}$, and
\item[(iv)] $V_0$ is invariant under the group of translation by integers.
\end{enumerate}

In view of property (iv), it is natural to ask the following question: Does there exist other groups of translations under which $V_0$ remains invariant? We shall consider the groups of translation by dyadic rationals. For $y\in\R$, Let $T_y$ be the (unitary) translation operator defined by $T_y f(x) = f(x-y)$. Consider the following groups of translation operators:
\[
\G_r =\{T_{\frac{m}{2^r}}:m\in\Z\},\ r\geq 0,\ r\in\Z, \quad{\rm and}\quad
\G_{\infty}=\{T_y:y\in\R\}.
\]
 
Let $\mathcal G$ be a set of bounded linear operators on $H^2(\R)$ and $V$ a closed subspace of $H^2(\R)$. We say that $V$ is $\mathcal G$-invariant if $Tf\in V$ for every $f\in V$ and $T\in{\mathcal G}$.
 
Let us denote by $\L_r$ the collection of all $H^2$-wavelets such that the corresponding space $V_0$ is $\G_r$-invariant. Clearly, we have the following inclusions:
\[
\L_{0}\supset\L_{1}\supset\L_{2}\supset\cdots\supset
\L_r\supset\L_{r+1}\supset\cdots\supset\L_{\infty}.
\]

We now define an equivalence relation on $H^2$-wavelets, where the equivalence classes are given by 
$\M_r = \L_r\setminus \L_{r+1}$, with $\M_\infty=\L_\infty$. Therefore, $\M_r$, $r\geq 0$ consists of $H^2$-wavelets for which $V_0$ is $\G_r$-invariant but not $\G_{r+1}$-invariant.

This equivalence relation was first defined in~\cite{Web} for the classical case of wavelets of $L^2(\R)$. In the same paper the author characterized the equivalence classes in terms of the support of the Fourier transform of the wavelets. He also proved that $\M_r$, $r=0,1,2,3$, are non-empty. Later, in~\cite{BM},~\cite{SW}, examples of wavelets of $L^2(\R)$ were constructed for each of these equivalence classes, by different methods.

The characterization of $\M_r$ can be easily carried over to the case of $H^2(\R)$. First of all we introduce some notation.

Let $\psi$ be an $H^2$-wavelet and let $E=\mbox{supp }\hat\psi$. For $k\in\Z$, define $E(k)=\{\xi\in E:\xi+2k\pi\in E\}$ and ${\mathcal E}_\psi=\{k\in\Z:E(k)\not=\emptyset\}$. Then the characterization of the equivalence classes is the following.

\begin{theorem}
\label{T.four}
\begin{itemize}
\item[(i)] $\M_\infty$ is precisely the collection of all $H^2$-MSF wavelets. 
\item[(ii)] An $H^2$-wavelet $\psi\in\M_r$, $r\geq 1$, if and only if  every element of    
 ${\mathcal E}_\psi$ is divisible by $2^r$ but there is an element of ${\mathcal E}_\psi$  
 not divisible by $2^{r+1}$.
\item[(iii)] An $H^2$-wavelet $\psi\in\M_0$ if and only if ${\mathcal E}_\psi$ contains an  odd integer.
\end{itemize}
\end{theorem}

The proof of the above theorem is an easy generalization of the corresponding result proved in~\cite{Web} for $L^2(\R)$. The purpose of this section is to show that all the equivalence classes are non-empty. Indeed, we show that, the $H^2$-wavelets constructed in the previous section serve as examples in $\M_r$, $r\geq 1$. To show that $\M_0$ is non-empty, we produce an interesting family of $H^2$-wavelet sets consisting of five disjoint intervals.

\begin{theorem}\label{T.five}
The equivalence classes $\M_r$, $r\in\N\cup\{0,\infty\}$, defined above, are non-empty.
\end{theorem}

\proof
We mentioned in the introduction about all previously known $H^2$-wavelets. Each of them is MSF. Hence, $\M_\infty$ is non-empty.

Now, fix $r\in\N$ and consider the $H^2$-wavelet $\psi_r$ defined in (\ref{E.four}). From Lemma~\ref{L.one}, we notice that ${\mathcal E}_{\psi_r}=\{0,\pm 2^r,\pm 2^{r+1}\}$. By Theorem~\ref{T.four}(ii), $\psi_r\in\M_r$.

We now construct an example in the equivalence class $\M_0$. For this purpose, it is natural to consider the case $r=0$ in (\ref{E.four}). Unfortunately this does not work since we get $\hat\psi_0=\chi_{[2\pi,4\pi]}$. Hence, $\psi_0$ is in $\M_\infty$, being an MSF wavelet.

Let $\pi<x<y<2\pi$ and $x+2\pi>2y$. That is, $(x,y)$ is in the interior of the triangle with vertices $(\pi, \frac{3}{2}\pi)$, $(\pi,2\pi)$ and $(2\pi, 2\pi)$. Consider the following set:
\[
K_{x,y}=[x,y]\cup[2\pi,2x]\cup[2y,x+2\pi]\cup[y+2\pi,4\pi]\cup[2x+4\pi,2y+4\pi].
\]
Let us denote the intervals in the right by $I_1,I_2,\dots,I_5$. The conditions on $x$ and $y$ ensure that these intervals are non-empty. Observe that $I_1$, $I_4-2\pi$, $I_2$, $I_5-4\pi$, $I_3$ are pairwise disjoint, and their union is $[x,x+2\pi]$. Similarly, 
$I_1$, $2^{-1}I_3$, $2^{-2}I_5$, $2^{-1}I_4$, $I_2$ are pairwise disjoint, and their union is $[x,2x]$. Hence, by Theorem~\ref{T.two}, $K_{x,y}$ is an $H^2$-wavelet set. 

In particular, we obtain a family of 5-interval $H^2$-wavelet sets where some of the endpoints of the intervals depend on two independent continuous real parameters.

Note that $2^{-1}I_3+2\pi$ is properly contained in $I_4$. Now, define the function $\psi_0$ by
\begin{eqnarray*}
\hat\psi_0(\xi)=
\begin{cases}
\frac{1}{\sqrt{2}} & {\rm if}~\xi\in I_3\cup (2^{-1}I_3)\cup(2^{-1}I_3+2\pi) \\
-\frac{1}{\sqrt 2} & {\rm if}~\xi\in I_3+4\pi \\
1 & {\rm if}~\xi\in K_{x,y}\setminus(2^{-1}I_3+2\pi) \\
0 & {\rm otherwise}.
\end{cases}
\end{eqnarray*}

It can be proved that $\psi_0$ is an $H^2$-wavelet. The proof is similar to that of Theorem~\ref{T.three} and we skip the proof in order to avoid repetition. It is also clear that ${\mathcal E}_{\psi_0}=\{0,\pm 1,\pm 2\}$. Hence by Theorem~\ref{T.four}(iii), $\psi_0\in\M_0$. This completes the proof.
\qed
\section{$H^2$-wavelets with Fourier transform discontinuous at the origin}

In this section we construct a family of wavelets for $H^2(\R)$ whose Fourier transforms are discontinuous at the origin. First we recall a result proved in~\cite{HW} for wavelets of $L^2(\R)$ (see Theorem~2.7, Chapter~3 of~\cite{HW}).
\begin{theorem}
Let $\psi$ be a wavelet for $L^2(\R)$ such that $\hat\psi$ has compact support and $|\hat\psi|$ is continuous at $0$. Then $\hat\psi=0$ a.e. in an open neighbourhood of the origin.
\end{theorem} 
This result also holds for $H^2(\R)$ with essentially the same proof. We are interested in the following question: Does there exist an $H^2$-wavelet such that $\hat\psi$ has compact support and does not vanish in any neighbourhood of the origin? In this section we shall give a positive answer to this question. We need the following concepts.

\begin{definition}
A set $A$ is said to be translation equivalent to a set $B$ if there exists a partition $\{A_n\}$ of $A$ and $k_n\in\Z$ such that $\{A_n+2k_n\pi\}$ is a partition of $B$. Similarly, $A$ is dilation equivalent to $B$ if there exists another partition $\{A_n'\}$ of $A$ and $j_n\in\Z$ such that $\{2^{j_n}A_n'\}$ is a partition of $B$.      
\end{definition}

Theorem \ref{T.two} has the following simple but useful consequence.

\begin{corollary}
\label{cor:wavset}
Let $K_1,K_2\subset\R^+$ and $K_1$ is both translation and dilation equivalent to $K_2$. Then $K_1$ is an $H^2$-wavelet set if and only if $K_2$ is so.
\end{corollary}

Let $r\in\N$ and $t_r=\tfrac{2^{r+1}}{2^{r+1}-1}\pi$. Then we know that
\[
K_r=[t_r,2\pi]\cup[2^{r+1}\pi, 2^{r+1}t_r]=I_r\cup J_r
\]
is an $H^2$-wavelet set (see equation (\ref{E.kr})). For $\epsilon>0$ such that 
$\epsilon < \tfrac{2^r-1}{2^{r+1}-1}\pi$, let
\begin{eqnarray*}
S_1 & = & \Bigl[\tfrac{t_r}{2}+\tfrac{\epsilon}{2^{r+1}},
          \tfrac{t_r}{2}+\epsilon\Bigr],\\
S_2 & = & [t_r+2\epsilon,2\pi],\\
S_3 & = & [2^{r+1}t_r, 2^{r+1}t_r+2\epsilon].
\end{eqnarray*}
The condition on $\epsilon$ ensures that $S_2$ a non-empty set. Let 
\[ 
\begin{array}{ll}
E_0 = S_1+2^{r+1}\pi, & 
F_0 = 2^{-(r+2)}E_0, \\
E_n = F_{n-1}+2^{r+1}\pi, & 
F_n = 2^{-(n+r+2)}E_n, \quad n\geq1,
\end{array} 
\]
Define the set
\begin{eqnarray}
\label{eqn:Wne}
K_{r,\epsilon} & = & 
\Bigl(J_r\setminus\bigcup_{n=0}^{\infty}E_n\Bigr)\cup
\Bigl(\bigcup_{n=0}^{\infty}F_n\Bigr)\cup (S_1\cup S_2 \cup S_3).
\end{eqnarray}

\begin{theorem}
For each $r\in\N$, the set $K_{r,\epsilon}$, defined in {\rm(\ref{eqn:Wne})}, 
is an $H^2$-wavelet set.
\end{theorem}

\proof
The result will follow from Corollary~\ref{cor:wavset} once we show that $K_{r,\epsilon}$ is both translation and dilation equivalent to the wavelet set $K_r$. First of all, we show by induction that $E_n\subset J_r$ for all $n\geq 0$.

Observe that $t_r+2^{r+1}\pi=2^{r+1}t_r$, hence, $[0, t_r]+2^{r+1}\pi=J_r$.
We have, $E_0= S_1+2^{r+1}\pi \subset[0,t_r]+2^{r+1}\pi=J_r$. 
Now assume that $E_m\subset J_r$. Then 
$F_m=2^{-(m+r+2)}E_m\subset 2^{-(m+1)}[0,t_r]\subset[0,t_r]$, 
hence, $E_{m+1} = F_m+2^{r+1}\pi\subset[0,t_r]+2^{r+1}\pi=J_r$.

The intervals $E_n$, $n\geq 0$, lie inside the interval $J_r$ and $E_{n+1}$ lies to the left of $E_n$ for all $n\geq 0$. Similarly, the intervals $F_n$, $n\geq 0$, lie in 
$2^{-(n+1)}[\pi,t_r]$ so that $F_{n+1}$ lies to the left of $F_n$ for $n\geq 0$. 

We now show that the set $K_{r,\epsilon}$ is dilation equivalent to the wavelet set $K_r$. We have
\begin{eqnarray*}
2S_1\cup S_2\cup\frac{1}{2^{r+1}}S_3 
& = &
  \Bigl[t_r+\frac{\epsilon}{2^r},t_r+2\epsilon\Bigr]
  \cup [t_r+2\epsilon,2\pi]
  \cup \Bigl[t_r,t_r+\frac{\epsilon}{2^r}\Bigr] \\
& = & [t_r,2\pi] = I_r,
\end{eqnarray*}
and 
\[
\Bigl(J_r\setminus\bigcup_{n=0}^{\infty}E_n\Bigr)       \cup\Bigl(\bigcup_{n=0}^{\infty}2^{n+r+2}F_n\Bigr)
= \Bigl(J_r\setminus\bigcup_{n=0}^{\infty}E_n\Bigr)\cup
\Bigl(\bigcup_{n=0}^{\infty}E_n\Bigr) = J_r,
\]
since $E_n\subset J_r$, for all $n\geq 0$, which proves the dilation equivalence.

Finally, we show that $K_{r,\epsilon}$ is translation equivalent to $K_r$. Observe that
\[
S_2\cup(S_3-2^{r+1}\pi) = [t_r+2\epsilon,2\pi]
\cup [t_r,t_r+2\epsilon] = I_r,
\]
and 
\begin{eqnarray*}
&   & 
\Bigl(J_r\setminus\bigcup_{n=0}^{\infty}E_n\Bigr)
\cup\Bigl(\bigcup_{n=0}^{\infty}(F_n+2^{r+1}\pi)\Bigr)
\cup(S_1+2^{r+1}\pi) \\
& = & 
\Bigl(J_r\setminus\bigcup_{n=0}^{\infty}E_n\Bigr)
\cup\Bigl(\bigcup_{n=1}^{\infty}E_n\Bigr)\cup E_0 = J_r.
\end{eqnarray*}
We have proved that $K_{r,\epsilon}$ is both dilation and translation equivalent to the $H^2$-wavelet set $K_r$. Therefore, by Corollary~\ref{cor:wavset}, it follows that $K_{r,\epsilon}$ is an $H^2$-wavelet set.
\qed

Let $\hat\psi_{r,\epsilon}$ be the characteristic function of the set $K_{r,\epsilon}$ 
so that $\psi_{r,\epsilon}$ is an $H^2$-wavelet. Since $F_n\subset2^{-(n+1)}[\pi,t_r]$ for all $n\geq 0$, $\hat\psi_{r,\epsilon}$ does not vanish in any neighbourhood of 0. In particular, it is discontinuous at the origin.


\begin{thebibliography}{amsplain}

\bibitem{ABM} Arcozzi, N., Behera, B., and Madan, S., 
{\it Large classes of minimally supported frequency wavelets of $L^2(\R)$ and $H^2(\R)$}. Preprint (2001).

\bibitem{AUS} Auscher, P.,
{\it Solution of two problems on wavelets},
J. Geom. Anal. {\bf 5} (1995), 181--236.

\bibitem{B} Behera, B., 
{\sl Band-limited wavelets and wavelet packets}, 
Ph. D. thesis, Indian Institute of Technology, Kanpur, India (2001).

\bibitem{BM} Behera, B., and Madan, S., 
{\it Characterization of a class of band-limited wavelets}, submitted.

\bibitem{DL} Dai, X., and Larson, D., 

\bibitem{HKLS} Ha, Y-H., Kang, H., Lee, J., and  Seo, J., 
{\it Unimodular wavelets for $L\sp 2$ and the Hardy space $H\sp 2$}, 
Michigan Math. J. {\bf 41} (1994), no. 2, 345--361. 

\bibitem{HW} Hern\'andez, E., and Weiss, G.,
{\sl A First Course on Wavelets}, CRC Press, Boca Raton (1996).
 
\bibitem{MAJ} Majchrowska, G., 
{\it Some new examples of wavelets in the Hardy space $H\sp 2({\Bbb R})$}, 
Bull. Polish Acad. Sci. Math. 49 (2001), no. 2, 141--149. 

\bibitem{SW} Schaffer, S., and  Weber, E., 
{\it Wavelets with the translation invariance property of order n}. Preprint (2001).

\bibitem{Web} Weber, E., 
{\it On the translation invariance of wavelet subspaces}, 
J. Fourier Anal. Appl. {\bf 6} (2000), no. 5, 551--558. 

\end{thebibliography}
\end{document}